\documentclass[11pt,twoside]{article}  
\usepackage[ansinew]{inputenc}
\usepackage{amsfonts}
\usepackage{amssymb,amsmath}
\usepackage{latexsym}

\newcommand {\debeq}	{\begin{eqnarray*}}
\newcommand {\fineq}	{\end{eqnarray*}}
\newcommand {\lbd}	{\lambda}
\newcommand     {\eps}  {\epsilon}
\newcommand     {\vareps}       {\varepsilon}

\newcommand	{\tendinfty}	
{\rightarrow\infty}

\newcommand	{\intgen}	
{\int_0^\infty}

\newcommand	{\PP}{\mathbb{P}}
\newcommand	{\EE}{\mathbb{E}}

\newtheorem	{thm}		{Theorem}[section]

\newtheorem	{lem} 	[thm]	{Lemma}

\newtheorem     {rem}           {Remark}
\newtheorem	{prop}	[thm]{Proposition}

\newcommand	{\indic}	[1]
{{\bf{1}}_{\{#1\}}}
\newcommand	{\indicbis}	[1]
{{\bf{1}}_{#1}}

\linethickness{0.4pt}

\begin{document}

\title{The allelic partition for coalescent point processes}
\author{\textsc{By Amaury Lambert
}
}
\date{\today}
\maketitle
\noindent\textsc{Laboratoire de Probabilités et Modèles Aléatoires\\
UMR 7599 CNRS and UPMC Univ Paris 06\\
Case courrier 188\\
4, Place Jussieu\\
F-75252 Paris Cedex 05, France}\\
\textsc{E-mail: }amaury.lambert@upmc.fr\\
\textsc{URL: }http://ecologie.snv.jussieu.fr/amaury/

\begin{abstract}
\noindent
Assume that individuals alive at time $t$ in some population can be ranked in such a way that the coalescence times between consecutive individuals are i.i.d. The ranked sequence of these branches is called a coalescent point process. We have shown in a previous work \cite{L} that splitting trees are important instances of such populations.

Here, individuals are given DNA sequences, and for a sample of $n$ DNA sequences belonging to distinct individuals, we consider the number $S_n$ of polymorphic sites (sites at which at least two sequences differ), and the number $A_n$ of distinct haplotypes (sequences differing at one site at least). 

It is standard to assume that mutations arrive at constant rate (on germ lines), and never hit the same site on the DNA sequence.
We study the mutation pattern associated with coalescent point processes under this assumption. Here, $S_n$ and $A_n$ grow linearly as $n$ grows, with explicit rate. However, when the branch lengths have infinite expectation, $S_n$ grows more rapidly, e.g. as $n \ln(n)$ for critical birth--death processes.

Then, we study the frequency spectrum of the sample, that is, the numbers of polymorphic sites/haplotypes carried by $k$ individuals in the sample. These numbers are shown to grow also linearly with sample size, and we provide simple explicit formulae for mutation frequencies and haplotype frequencies. For critical birth--death processes, mutation frequencies are given by the harmonic series and haplotype frequencies by Fisher's logarithmic series.

\end{abstract}  	
\medskip
\textit{Running head.} The allelic partition for coalescent point processes.

\section{Introduction}

\subsection{The coalescent point process}

Splitting trees are those random trees where individuals give birth at constant rate $b$ during a lifetime with general distribution $\Lambda(\cdot)/b$, to i.i.d. copies of themselves (see \cite{GK}), where $\Lambda$ is a positive measure on $(0,\infty]$ with total mass $b$ called the \emph{lifespan measure}. In \cite{L}, we have shown that if the splitting tree is started from one individual with known birth time, say $0$, and known death time, then individuals alive at time $t$ can be ranked in such a way that the coalescence times between consecutive individuals are i.i.d.

\begin{figure}[ht]

\unitlength 2mm 
\linethickness{0.4pt}
\begin{picture}(66,33)(-5,10)
\put(4,39.875){\line(1,0){62}}
\put(10,40){\line(0,-1){9}}
\put(14,40){\line(0,-1){11.5}}
\put(18,40){\line(0,-1){4}}
\put(22,40){\line(0,-1){7}}
\put(26,40){\line(0,-1){16}}
\put(30,40){\line(0,-1){6}}
\put(34,39.875){\line(0,-1){8.5}}
\put(38,40){\line(0,-1){5.5}}
\put(42,40){\line(0,-1){11.5}}
\put(46,40){\line(0,-1){3.625}}
\put(50,40){\line(0,-1){22}}
\put(54,40){\line(0,-1){5}}
\put(58,40){\line(0,-1){7.5}}
\put(62,39.875){\line(0,-1){5}}
\put(6,40){\line(0,-1){25}}
\put(5.93,14.93){\line(0,-1){.8}}
\put(5.93,13.33){\line(0,-1){.8}}
\put(5.93,11.73){\line(0,-1){.8}}
\put(9.93,30.93){\line(-1,0){.8}}
\put(8.33,30.93){\line(-1,0){.8}}
\put(6.73,30.93){\line(-1,0){.8}}
\put(13.93,28.43){\line(-1,0){.8889}}
\put(12.152,28.43){\line(-1,0){.8889}}
\put(10.374,28.43){\line(-1,0){.8889}}
\put(8.596,28.43){\line(-1,0){.8889}}
\put(6.819,28.43){\line(-1,0){.8889}}
\put(17.805,35.93){\line(-1,0){.8}}
\put(16.205,35.93){\line(-1,0){.8}}
\put(14.605,35.93){\line(-1,0){.8}}
\put(21.93,32.93){\line(-1,0){.8889}}
\put(20.152,32.93){\line(-1,0){.8889}}
\put(18.374,32.93){\line(-1,0){.8889}}
\put(16.596,32.93){\line(-1,0){.8889}}
\put(14.819,32.93){\line(-1,0){.8889}}
\put(25.93,23.93){\line(-1,0){.9524}}
\put(24.025,23.93){\line(-1,0){.9524}}
\put(22.12,23.93){\line(-1,0){.9524}}
\put(20.215,23.93){\line(-1,0){.9524}}
\put(18.311,23.93){\line(-1,0){.9524}}
\put(16.406,23.93){\line(-1,0){.9524}}
\put(14.501,23.93){\line(-1,0){.9524}}
\put(12.596,23.93){\line(-1,0){.9524}}
\put(10.692,23.93){\line(-1,0){.9524}}
\put(8.787,23.93){\line(-1,0){.9524}}
\put(6.882,23.93){\line(-1,0){.9524}}
\put(29.93,33.93){\line(-1,0){.8}}
\put(28.33,33.93){\line(-1,0){.8}}
\put(26.73,33.93){\line(-1,0){.8}}
\put(33.93,31.43){\line(-1,0){.8889}}
\put(32.152,31.43){\line(-1,0){.8889}}
\put(30.374,31.43){\line(-1,0){.8889}}
\put(28.596,31.43){\line(-1,0){.8889}}
\put(26.819,31.43){\line(-1,0){.8889}}
\put(37.93,34.43){\line(-1,0){.8}}
\put(36.33,34.43){\line(-1,0){.8}}
\put(34.73,34.43){\line(-1,0){.8}}
\put(41.93,28.43){\line(-1,0){.9412}}
\put(40.047,28.43){\line(-1,0){.9412}}
\put(38.165,28.43){\line(-1,0){.9412}}
\put(36.283,28.43){\line(-1,0){.9412}}
\put(34.4,28.43){\line(-1,0){.9412}}
\put(32.518,28.43){\line(-1,0){.9412}}
\put(30.636,28.43){\line(-1,0){.9412}}
\put(28.753,28.43){\line(-1,0){.9412}}
\put(26.871,28.43){\line(-1,0){.9412}}
\put(45.93,36.43){\line(-1,0){.8}}
\put(44.33,36.43){\line(-1,0){.8}}
\put(42.73,36.43){\line(-1,0){.8}}
\put(49.93,17.93){\line(-1,0){.9778}}
\put(47.974,17.93){\line(-1,0){.9778}}
\put(46.019,17.93){\line(-1,0){.9778}}
\put(44.063,17.93){\line(-1,0){.9778}}
\put(42.107,17.93){\line(-1,0){.9778}}
\put(40.152,17.93){\line(-1,0){.9778}}
\put(38.196,17.93){\line(-1,0){.9778}}
\put(36.241,17.93){\line(-1,0){.9778}}
\put(34.285,17.93){\line(-1,0){.9778}}
\put(32.33,17.93){\line(-1,0){.9778}}
\put(30.374,17.93){\line(-1,0){.9778}}
\put(28.419,17.93){\line(-1,0){.9778}}
\put(26.463,17.93){\line(-1,0){.9778}}
\put(24.507,17.93){\line(-1,0){.9778}}
\put(22.552,17.93){\line(-1,0){.9778}}
\put(20.596,17.93){\line(-1,0){.9778}}
\put(18.641,17.93){\line(-1,0){.9778}}
\put(16.685,17.93){\line(-1,0){.9778}}
\put(14.73,17.93){\line(-1,0){.9778}}
\put(12.774,17.93){\line(-1,0){.9778}}
\put(10.819,17.93){\line(-1,0){.9778}}
\put(8.863,17.93){\line(-1,0){.9778}}
\put(6.907,17.93){\line(-1,0){.9778}}
\put(53.93,34.93){\line(-1,0){.8}}
\put(52.33,34.93){\line(-1,0){.8}}
\put(50.73,34.93){\line(-1,0){.8}}
\put(57.93,32.43){\line(-1,0){.8889}}
\put(56.152,32.43){\line(-1,0){.8889}}
\put(54.374,32.43){\line(-1,0){.8889}}
\put(52.596,32.43){\line(-1,0){.8889}}
\put(50.819,32.43){\line(-1,0){.8889}}
\put(61.93,34.93){\line(-1,0){.8}}
\put(60.33,34.93){\line(-1,0){.8}}
\put(58.73,34.93){\line(-1,0){.8}}
\put(6,41){\makebox(0,0)[cc]{$0$}}
\put(10,41){\makebox(0,0)[cc]{$1$}}
\put(14,41){\makebox(0,0)[cc]{$2$}}
\put(18,41){\makebox(0,0)[cc]{$3$}}
\put(22,41){\makebox(0,0)[cc]{$4$}}
\put(26,41){\makebox(0,0)[cc]{$5$}}
\put(30,41){\makebox(0,0)[cc]{$6$}}
\put(34,41){\makebox(0,0)[cc]{$7$}}
\put(38,41){\makebox(0,0)[cc]{$8$}}
\put(42,41){\makebox(0,0)[cc]{$9$}}
\put(46,41){\makebox(0,0)[cc]{$10$}}
\put(50,41){\makebox(0,0)[cc]{$12$}}
\put(54,41){\makebox(0,0)[cc]{$13$}}
\put(58,41){\makebox(0,0)[cc]{$14$}}
\put(62,41){\makebox(0,0)[cc]{$15$}}
\end{picture}

\caption{ A coalescent point process for $n=16$ individuals.}
\label{fig : coalpointproc}
\end{figure}
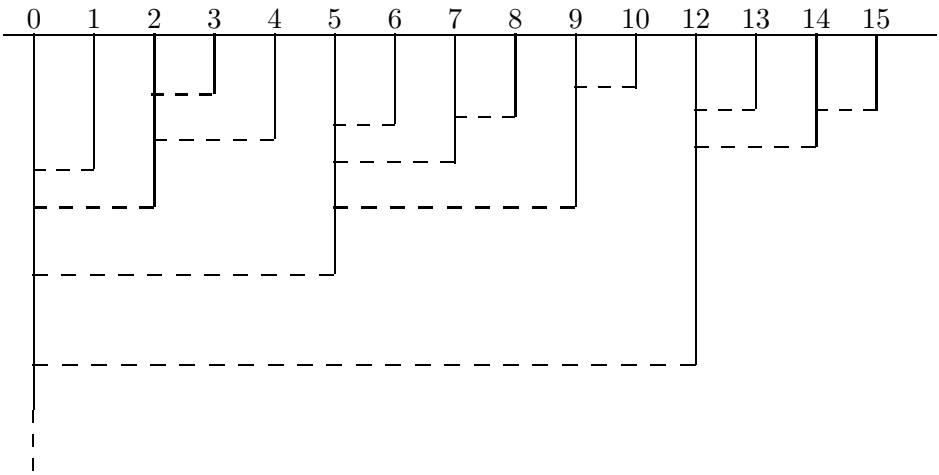

Specifically, let $N_t$ be the number of individuals alive at time $t$. The process $(N_t;t\ge 0)$ is a (homogeneous, binary) Crump--Mode--Jagers process, and is not Markovian unless $\Lambda$ has an exponential density or is a point mass at $\infty$. To these $N_t$ individuals, give labels $0,1,\ldots, N_t-1$ according to the (unique) order complying with the following rule : `any individual comes before her own descendants, but after her younger siblings and their descendants'. For any integers $i,k$ such that $0\le i<i+k< N_t$, we let $C_{i,i+k}$ be the \emph{coalescence time} (or \emph{divergence time}) between individual $i$ and individual $i+k$, that is, the time elapsed since the lineages of individual $i$ and $i+k$ have diverged. Also define $H_{i+1}:=C_{i,i+1}$. Then recall from \cite{L} that for a splitting tree,
\begin{equation}
\label{eqn : def coal}
C_{i,i+k}=\max\{H_{i+1},\ldots,H_{i+k}\} 
\end{equation}
and conditional on $\{N_t\not=0\}$, the sequence $(H_i;1\le i \le N_t-1)$ has the same law as a sequence of i.i.d. r.v. killed at its first value  $\ge t$. As a by-product, we get that the law of $N_t$ conditional on $\{N_t\not=0\}$ is geometric.

The aforementioned property comes from the fact that the jumping contour process of the splitting tree is a Lévy process $X=(X_s;s\ge0)$ with Lévy measure $\Lambda$ and drift coefficient $-1$. Then the excursions of the contour process between consecutive visits of points at height $t$ are i.i.d. excursions of $X$. As a consequence, the $(H_i)$ are also i.i.d., and their common distribution is that of $H':=t-\inf_s X_s$, where $X$ is started at $t$ and killed upon hitting $\{0\}\cup(t,+\infty)$. Note that all branch lengths but the last one are distributed as some r.v. $H$ which is $H'$ conditioned to be smaller than $t$.
The distribution of $H'$ can be expressed in terms of a nonnegative, nondecreasing, differentiable function $W$, called the \emph{scale function} of $X$, such that $W(0)=1$
\begin{equation}
\label{eqn : scale}
\PP(H'> x)=\frac{1}{W(x)}\qquad x\ge 0.
\end{equation}
The scale function $W$ is characterised by its Laplace transform (see e.g. \cite{B})
\begin{equation}
\label{eqn : LT scale}
\intgen dx\, e^{-\lbd x} \, W(x) = \left(\lbd -\intgen \Lambda(dx) (1-e^{-\lbd x}) \right)^{-1}.
\end{equation}
From now on, with no need to refer to the framework of splitting trees, we will consider the genealogy of what we call a \emph{coalescent point process} (originating from \cite{P} where $\Lambda (dx)=b^2\exp(-bx) dx$) :
\begin{enumerate}
\item let $H_1, H_2,\ldots$ be a sequence of independent random variables called \emph{branch lengths} all distributed as some positive r.v. $H$, and set $H_0$ to equal $+\infty$.
\item
the genealogy of the population $\{0,1,2,\ldots\}$ is given by \eqref{eqn : def coal}.
\end{enumerate}
We will stick to the notation
$$
W(x):=\frac{1}{\PP(H> x)}\qquad x\ge 0.
$$
It will always be implicit that a \emph{sample} of $n$ individuals refers to the  \emph{first} $n$ individuals $\{0,1,\ldots,n-1\}$.
\begin{rem}
\label{rem : H' et H}
In the case of splitting trees, conditional on $\{N_t\not=0\}$, $N_t$ is geometric with success probability $\PP(H'>t)$, and conditional on $\{N_t=n\}$, the branch lengths $(H_i;1\le i \le n-1)$ are i.i.d. with distribution $\PP(H'\in\cdot\mid H'<t)$. In what follows, we will repeatedly refer to the genealogy of a splitting tree with $n$ leaves by setting the r.v. $H$ to equal $H'$, without the conditioning (i.e. $t\to\infty$). In the subcritical case, this amounts to considering quasi-stationary populations, which are those populations conditioned to be still alive at time $t$, as $t\tendinfty$ (see e.g. \cite{L9}). Another possibility would be, as in \cite{AP}, to give a prior distribution to the time $t$ of origin, and condition the whole tree on $\{N_t=n\}$. Then as $n\to\infty$, the posterior distribution of $t$ goes to $\infty$, and we would be left with a (possibly different) distribution of $H$ charging the whole half-line. 
\end{rem}
\begin{rem}
No distribution of edge lengths can make the coalescent point process coincide with the Kingman coalescent \cite{Ki}. Indeed, here, the smallest branch length in a sample of $n$ individuals is the minimum of  $n-1$ i.i.d. random variables, whereas in the Kingman coalescent, it is the minimum of $n(n-1)/2$ i.i.d. random variables (with exponential distribution).
\end{rem}

Our goal is to characterise the mutation pattern for samples of $n$ individuals, mainly as $n$ gets large. We specify the mutation scheme in the next subsection. 
 
Works studying mutation patterns arising from random genealogies are numerous. Mutation patterns related to populations with fixed size (Wright--Fisher model, Kingman coalescent) are well-known and culminate in \emph{Ewens' sampling formula} (see \cite{D1} for a comprehensive account on that subject). More recent works concern mutation patterns related to more general coalescents \cite{BG,M}, to branching populations \cite{AD, Bclusters}, or to both \cite{BBS}.

\subsection{Mutation scheme}
We adopt two classical assumptions on mutation schemes from population genetics (see e.g. \cite{Ew})
\begin{enumerate}
\item 
mutations occur at \emph{constant rate} $\theta$ on germ lines,
\item
mutations are \emph{neutral}, that is, they have no effect on birth rates and lifetimes.
\end{enumerate}
As is usual, we assume that mutations are point substitutions occurring at a single site on the DNA sequence, and that each site can be hit at most once by a mutation. This last assumption is known as the \emph{infinitely-many sites model} (ISM). Instances of DNA sequence are called \emph{alleles} or \emph{haplotypes}, so that under the ISM, each mutation yields a new allele. Without reference to DNA sequences, this last assumption by itself is known as the \emph{infinitely-many alleles model} (IAM).\\
\\
Specifically, we let $({\cal P}_i;i=0,1,2\ldots)$ be independent Poisson measures on $(0,\infty)$ with intensity $\theta$ (cf. assumption 1). For each $i$ we denote the atoms of ${\cal P}_i$ by $\ell_{i1}<\ell_{i2}<\cdots$ and call them \emph{mutations}. Now let $H_1, H_2, \ldots$ be an independent coalescent point process (cf. assumption 2). In agreement with the genealogical structure of a coalescent point process explained in the beginning of this section, we will say that individual $i+k$ \emph{carries} (or \emph{bears}) mutation $\ell_{ij}$ if $k\ge 0$ and
$$
\max\{H_{i+1},\ldots H_{i+k} \}< \ell_{ij} <H_i,
$$
where we agree that $\max\emptyset=0$ and $H_0=+\infty$. The second inequality is trivially due to the fact that we throw away all atoms $\ell_{ij}$ such that $H_i\le\ell_{ij}$. The set of mutations that an individual bears is her \emph{allele} or her \emph{haplotype}, or merely her type.\\
\\
For a sample of $n$ individuals, we call $S_n$ the number of \emph{polymorphic sites}, that is, the number of mutations $(\ell_{ij}; 0\le i\le n-1, j\ge 1)$ that are carried by at least one individual and at most $n-1$. Formally, this yields
$$
S_n=\mbox{Card}\{\ell_{ij}< H_i, 1\le i\le n-1, j\ge 1\}+\mbox{Card}\{\ell_{0j}<\max\{H_1,\ldots, H_{n-1}\}, j\ge 1\}.
$$
Further, we define $S_n(k)$ as the number of mutations carried by $k$ individuals in the sample. In particular,
$$
S_n=\sum_{k=1}^{n-1}S_n(k).
$$
The sequence $(S_n(1),\ldots, S_n(n-1))$ is called the \emph{site frequency spectrum} of the sample.

\begin{figure}[ht]

\unitlength 3mm 
\linethickness{0.4pt}
\begin{picture}(61,33)(3,3)
\put(6,30){\line(1,0){45}}
\put(12,30){\line(0,-1){6}}
\put(22,30){\line(0,-1){6.5}}
\put(27,30){\line(0,-1){3.5}}
\put(32,30){\line(0,-1){7.5}}
\put(42,30){\line(0,-1){6}}
\put(47,30){\line(0,-1){3}}
\put(7,30){\line(0,-1){24}}
\put(17,30){\line(0,-1){12}}
\put(37,30){\line(0,-1){16.5}}
\put(7,15.5){\circle*{.707}}
\put(17,19.5){\circle*{.707}}
\put(22,28){\circle*{.707}}
\put(37,22){\circle*{.707}}
\put(37,17){\circle*{.707}}
\put(42,25){\circle*{.707}}
\put(11.93,23.93){\line(-1,0){.8333}}
\put(10.263,23.93){\line(-1,0){.8333}}
\put(8.596,23.93){\line(-1,0){.8333}}
\put(16.93,17.93){\line(-1,0){.9091}}
\put(15.112,17.93){\line(-1,0){.9091}}
\put(13.293,17.93){\line(-1,0){.9091}}
\put(11.475,17.93){\line(-1,0){.9091}}
\put(9.657,17.93){\line(-1,0){.9091}}
\put(7.839,17.93){\line(-1,0){.9091}}
\put(21.93,23.43){\line(-1,0){.8438}}
\put(20.242,23.43){\line(-1,0){.8438}}
\put(18.555,23.43){\line(-1,0){.8438}}
\put(26.93,26.43){\line(-1,0){.8333}}
\put(25.263,26.43){\line(-1,0){.8333}}
\put(23.596,26.43){\line(-1,0){.8333}}
\put(31.93,22.43){\line(-1,0){.9375}}
\put(30.055,22.43){\line(-1,0){.9375}}
\put(28.18,22.43){\line(-1,0){.9375}}
\put(26.305,22.43){\line(-1,0){.9375}}
\put(24.43,22.43){\line(-1,0){.9375}}
\put(22.555,22.43){\line(-1,0){.9375}}
\put(20.68,22.43){\line(-1,0){.9375}}
\put(18.805,22.43){\line(-1,0){.9375}}
\put(36.93,13.43){\line(-1,0){.9677}}
\put(34.994,13.43){\line(-1,0){.9677}}
\put(33.059,13.43){\line(-1,0){.9677}}
\put(31.123,13.43){\line(-1,0){.9677}}
\put(29.188,13.43){\line(-1,0){.9677}}
\put(27.252,13.43){\line(-1,0){.9677}}
\put(25.317,13.43){\line(-1,0){.9677}}
\put(23.381,13.43){\line(-1,0){.9677}}
\put(21.446,13.43){\line(-1,0){.9677}}
\put(19.51,13.43){\line(-1,0){.9677}}
\put(17.575,13.43){\line(-1,0){.9677}}
\put(15.639,13.43){\line(-1,0){.9677}}
\put(13.704,13.43){\line(-1,0){.9677}}
\put(11.768,13.43){\line(-1,0){.9677}}
\put(9.833,13.43){\line(-1,0){.9677}}
\put(7.897,13.43){\line(-1,0){.9677}}
\put(41.93,23.93){\line(-1,0){.8333}}
\put(40.263,23.93){\line(-1,0){.8333}}
\put(38.596,23.93){\line(-1,0){.8333}}
\put(46.93,26.93){\line(-1,0){.8333}}
\put(45.263,26.93){\line(-1,0){.8333}}
\put(43.596,26.93){\line(-1,0){.8333}}
\put(7,8){\circle*{.707}}
\put(8,15.5){\makebox(0,0)[cc]{$c$}}
\put(8,8){\makebox(0,0)[cc]{$g$}}
\put(18,19.5){\makebox(0,0)[cc]{$h$}}
\put(23,28){\makebox(0,0)[cc]{$e$}}
\put(38,22){\makebox(0,0)[cc]{$f$}}
\put(38,17){\makebox(0,0)[cc]{$b$}}
\put(43,25){\makebox(0,0)[cc]{$d$}}
\put(7,31){\makebox(0,0)[cc]{$gc$}}
\put(12,31){\makebox(0,0)[cc]{$gc$}}
\put(17,31){\makebox(0,0)[cc]{$gch$}}
\put(22,31){\makebox(0,0)[cc]{$gche$}}
\put(27,31){\makebox(0,0)[cc]{$gch$}}
\put(32,31){\makebox(0,0)[cc]{$gch$}}
\put(37,31){\makebox(0,0)[cc]{$gbf$}}
\put(42,31){\makebox(0,0)[cc]{$gbfd$}}
\put(47,31){\makebox(0,0)[cc]{$gbfd$}}
\put(7,33){\makebox(0,0)[cc]{$0$}}
\put(12,33){\makebox(0,0)[cc]{$1$}}
\put(17,33){\makebox(0,0)[cc]{$2$}}
\put(22,33){\makebox(0,0)[cc]{$3$}}
\put(27,33){\makebox(0,0)[cc]{$4$}}
\put(32,33){\makebox(0,0)[cc]{$5$}}
\put(37,33){\makebox(0,0)[cc]{$6$}}
\put(42,33){\makebox(0,0)[cc]{$7$}}
\put(47,33){\makebox(0,0)[cc]{$8$}}
\put(33,8.5){\framebox(1.5,1.5)[cc]{}}
\put(33,4){\framebox(1.5,1.5)[cc]{}}
\put(35,8.5){\framebox(1.5,1.5)[cc]{}}
\put(35,4){\framebox(1.5,1.5)[cc]{}}
\put(37,8.5){\framebox(1.5,1.5)[cc]{}}
\put(37,4){\framebox(1.5,1.5)[cc]{}}
\put(39,8.5){\framebox(1.5,1.5)[cc]{}}
\put(39,4){\framebox(1.5,1.5)[cc]{}}
\put(41,8.5){\framebox(1.5,1.5)[cc]{}}
\put(41,4){\framebox(1.5,1.5)[cc]{}}
\put(43,8.5){\framebox(1.5,1.5)[cc]{}}
\put(43,4){\framebox(1.5,1.5)[cc]{}}
\put(45,8.5){\framebox(1.5,1.5)[cc]{}}
\put(45,4){\framebox(1.5,1.5)[cc]{}}
\put(47,8.5){\framebox(1.5,1.5)[cc]{}}
\put(47,4){\framebox(1.5,1.5)[cc]{}}
\put(33.75,7.25){\makebox(0,0)[cb]{$a$}}
\put(33.75,2.75){\makebox(0,0)[cb]{$a$}}
\put(35.75,7.25){\makebox(0,0)[cb]{$b$}}
\put(35.75,2.75){\makebox(0,0)[cb]{$b$}}
\put(37.75,7.25){\makebox(0,0)[cb]{$c$}}
\put(37.75,2.75){\makebox(0,0)[cb]{$c$}}
\put(39.75,7.25){\makebox(0,0)[cb]{$d$}}
\put(39.75,2.75){\makebox(0,0)[cb]{$d$}}
\put(41.75,7.25){\makebox(0,0)[cb]{$e$}}
\put(41.75,2.75){\makebox(0,0)[cb]{$e$}}
\put(43.75,7.25){\makebox(0,0)[cc]{$f$}}
\put(43.75,2.75){\makebox(0,0)[cc]{$f$}}
\put(45.75,7.25){\makebox(0,0)[cc]{$g$}}
\put(45.75,2.75){\makebox(0,0)[cc]{$g$}}
\put(47.75,7.25){\makebox(0,0)[cb]{$h$}}
\put(47.75,2.75){\makebox(0,0)[cb]{$h$}}
\put(31.5,9){\makebox(0,0)[rc]{Haplotype of individual 3}}
\put(31.5,4.5){\makebox(0,0)[rc]{Haplotype of individual 6}}
\multiput(45,10)(.03333333,-.03333333){45}{\line(0,-1){.03333333}}
\multiput(41,10)(.03333333,-.03333333){45}{\line(0,-1){.03333333}}
\multiput(37,10)(.03333333,-.03333333){45}{\line(0,-1){.03333333}}
\multiput(47,10)(.03333333,-.03333333){45}{\line(0,-1){.03333333}}
\multiput(45,5.5)(.03333333,-.03333333){45}{\line(0,-1){.03333333}}
\multiput(43,5.5)(.03333333,-.03333333){45}{\line(0,-1){.03333333}}
\multiput(35,5.5)(.03333333,-.03333333){45}{\line(0,-1){.03333333}}
\multiput(45,8.5)(.03333333,.03333333){45}{\line(0,1){.03333333}}
\multiput(41,8.5)(.03333333,.03333333){45}{\line(0,1){.03333333}}
\multiput(37,8.5)(.03333333,.03333333){45}{\line(0,1){.03333333}}
\multiput(47,8.5)(.03333333,.03333333){45}{\line(0,1){.03333333}}
\multiput(45,4)(.03333333,.03333333){45}{\line(0,1){.03333333}}
\multiput(43,4)(.03333333,.03333333){45}{\line(0,1){.03333333}}
\multiput(35,4)(.03333333,.03333333){45}{\line(0,1){.03333333}}
\multiput(48.93,9.18)(.8,0){6}{{\rule{.4pt}{.4pt}}}
\multiput(48.93,4.742)(.8,0){6}{{\rule{.4pt}{.4pt}}}
\end{picture}
\caption{ A coalescent point process with mutations for a sample of $n=9$ individuals. Site $a$ is \emph{not} polymorphic because \emph{no} individual in the sample carries a mutation at that site; site $g$ is \emph{not} polymorphic because \emph{all} individuals in the sample carry the mutation at that site. The number of polymorphic sites is $S_n=6$. The number of distinct haplotypes is $A_n=5$.}
\label{fig : mutations}
\end{figure}
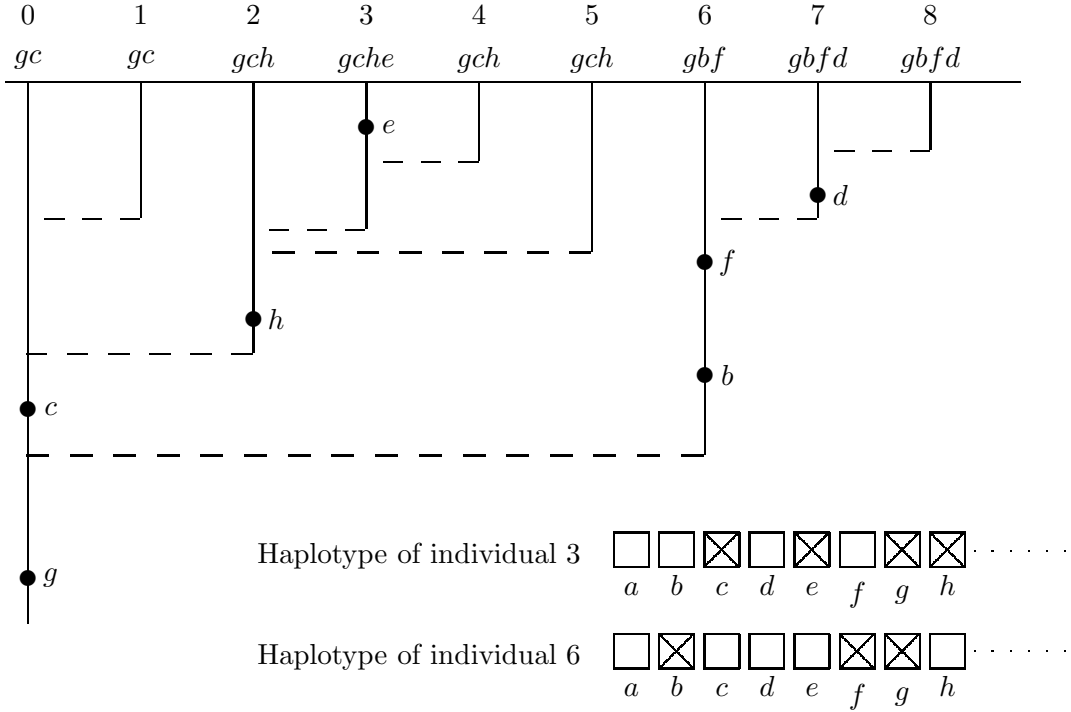

Similarly, we call $A_n$ the number of \emph{distinct haplotypes} in a sample of $n$ individuals, that is, the number of alleles that are carried by at least one individual, and $A_n(k)$ as the number of alleles carried by $k$ individuals. In particular, we have
$$
A_n=\sum_{k=1}^{n}A_n(k) \mbox{ and }\sum_{k=1}^n kA_n(k) = n.
$$
The sequence $(A_n(1),\ldots, A_n(n))$ is called the \emph{allele frequency spectrum} of the sample. 

\begin{rem}
One always has the inequality $S_n\ge A_n -1$. Indeed, apart from the ancestral haplotype, each new haplotype independent of at least one new mutation.
\end{rem}

\subsection{Examples of coalescent point processes}
Before going into the main part of this work, we provide a few simple examples of coalescent point processes derived from splitting trees, in part for application purposes.

\paragraph{Yule tree.} When $\Lambda$ is a point mass at $\infty$, the splitting tree is a Yule tree, and $(N_t;t\ge0)$ is a pure-birth binary process with birth rate, say $a$. Then $W(x)=e^{ax}$, and $H$ has an exponential distribution with parameter $a$ (see \cite{P}). 

\paragraph{Birth--death process.}
When $\Lambda$ has an exponential density, $(N_t;t\ge 0)$ is a Markovian birth--death process with (birth rate $b$ and) death rate, say $d$. Then it is known (see \cite{L} for example) that if $b\not=d$, then
$$
W(x) = \frac{d-be^{(b-d)x}}{d-b}\qquad x\ge 0,
$$
whereas if $b=d=:a$,
$$
W(x)=1+ax\qquad x\ge 0.
$$
Notice that in the subcritical case ($b<d$), $H$ can take the value $\infty$ with probability $1-(b/d)$, which is due to the constrained size of quasi-stationary populations (see Remark \ref{rem : H' et H}). Elementary calculations show that $H$ conditioned to be finite has the same law as the branch length of a \emph{supercritical} birth--death process with birth rate $d$ and death rate $b$.

\paragraph{Consistency and sampling.}
The genealogy associated with a coalescent point process is \emph{consistent} in the sense that the genealogy of a sample of $n$ individuals has the same law as that of a sample of $n+1$ individuals from which the \emph{last} individual has been withdrawn (in the splitting tree framework, the last individual is the individual who has no descendants in the sample, and whose ancestors have no elder sibling with descendants in the sample). This property would not hold any longer if the withdrawn individual was chosen at random. 

On the other hand, if all individuals in the population are censused \emph{independently} with probability $c$, then the genealogy of the census is still that of a coalescent point process. Indeed, the typical branch length is $H''$, where
$$
H''\stackrel{\cal L}{=}\max\{H_1,\ldots, H_K\},
$$
and $K$ is an independent (modified) geometric r.v., that is, $\PP(K=j)=c(1-c)^{j-1}$. As a consequence,
$$
\frac{1}{W_c(x)}:=\PP(H''> x)=1-\sum_{j\ge 1}c(1-c)^{j-1}\PP(H\le x)^{j}\qquad x\ge 0.
$$
This last equation also reads
$$
W_c = 1-c +cW.
$$
Applying this Bernoulli sampling procedure with intensity $c$ to the previous examples yields the following elementary results.
\begin{description}
	\item[--] the census of a Yule population has the genealogy of a birth--death process population, with birth rate $ac$ and death rate $a(1-c)$
	\item[--] the census of a birth--death process population has the genealogy of another birth--death process population with birth rate $bc$ and death rate $d-b(1-c)$. In particular, censusing a critical birth--death process population with rate $b=d=:a$ amounts to replacing $a$ with $ac$. 
\end{description}

\paragraph{Infinite lifespan measure.}
Actually, everything that was stated about splitting trees still holds if the lifespan measure is infinite, provided the lifespans of children remain summable, that is $\intgen(1\wedge r)\Lambda(dr)<\infty$. In particular, one still has $W(0)=1$, and the number of individuals alive at a fixed time remains a.s. finite.

On the contrary, it is a completely different task to define the real tree whose jumping contour process is a Lévy process with no negative jumps but \emph{infinite variation} (see \cite{B}). However, in our setting, this only requires replacing the coalescent point process $H_1, H_2,\ldots$ with a true Poisson point process with intensity measure $ds\,\nu(dx)$, where $\nu$ is a $\sigma$-finite positive measure defined as the push forward of the excursion measure of $X$ away from $\{t\}$ by the function which maps an excursion $\eps$ into $t-\inf_s\eps_s$. Similarly as in the finite variation case, 
$$
W(x):=\frac{1}{\nu((x,\infty))}\qquad x\ge0.
$$
In the Brownian case, for example $\nu(dx)= x^{-2}dx$ (again, see \cite{P}), that is, $W(x)=x$.

Here, the analogue of Bernoulli sampling with intensity $c$ consists in taking the maximum $H''$ of the point process on an interval with exponential length of parameter $c$ (instead of a geometric length). Now $c$ can take any positive value. Standard calculations then yield
$$
\frac{1}{W_c(x)}:=\PP(H''> x)=1-\frac{c}{c+\nu((x,\infty))}\qquad x\ge 0,
$$
so that
$$
W_c=1+cW.
$$
As far as splitting trees with infinite variation are concerned, we will only focus on the stable case, where $W(x)=x^{\alpha-1}$ for some $\alpha\in(1,2]$, the Brownian case corresponding to $\alpha=2$. In particular, we see that the Brownian coalescent point process censused with intensity $c$ has the same law as the coalescent point process associated with a critical birth--death process with rate $c$.

\subsection{Statements, outline, examples}
Our results regarding polymorphic sites are stated in Section \ref{sec : S}. 

In the first two subsections of Section \ref{sec : S}, we assume that $\EE(H)$ is finite. Theorem \ref{thm : LLN+CLT} provides a law of large numbers and a central limit theorem  (if $H$ has a second moment) on the number $S_n$ of polymorphic sites. In particular, 
\begin{equation}
\lim_{n\tendinfty}\frac{S_n}{n}=\theta\,\EE(H)\qquad\mbox{ a.s. }
\end{equation}
We also give exact explicit formulae for the \emph{expectation} of the number $S_n(k)$ of mutations carried by $k$ individuals in a sample of $n$.

In the third subsection, we make the less stringent assumption that $\EE(\min (H_1, H_2))$ is finite. Theorem \ref{thm : SFS large} then gives the asymptotic behaviour of the site frequency spectrum of large samples via the following a.s. convergence 
\begin{equation}
\label{eqn : SFS large}
\lim_{n\tendinfty}\frac{S_n(k)}{n} =	\theta\,\intgen \frac{dx}{W(x)^2}\left(1-\frac{1}{W(x)}\right)^{k-1}.
\end{equation}

In the fourth subsection, we treat the case of stable laws with parameter $\alpha$, that is, $W$ is given by $W(x)=1+c x^{\alpha -1}$, where $\alpha\in(1,2]$ and $c$ is some positive parameter that can be interpreted as a sampling intensity. Since here $\EE(H)=\infty$, the only result holding in the stable case is \eqref{eqn : SFS large}, and only for $\alpha>3/2$. Theorems \ref{thm : brown} and \ref{thm : stable} give the asymptotic behaviour of $S_n$. When $\alpha=2$, $S_n/n\ln (n)$ converges in probability (to $\theta/c$), and when $\alpha\not=2$, $S_n/n^{\beta}$ converges in distribution, with $\beta=1/(\alpha-1)$.\\
\\
Section \ref{sec : A} displays our results regarding distinct haplotypes. The trick is to characterise the law of the branch length $H^\theta$ of the next individual bearing no mutation other than those carried by, say, individual 0. Proposition \ref{prop : next branch with no} does this as follows
$$
\frac{1}{\PP(H_\theta>x)}=:W_\theta(x)=1+\int_0^x W'(u)e^{-\theta u}\, du\qquad x\ge 0.
$$
Theorem \ref{thm : AFS large} states a.s. convergences without moment existence assumptions. Specifically,
\begin{equation}
\lim_{n\tendinfty}\frac{A_n}{n}=\EE\left(1-e^{-\theta H^\theta}\right)\qquad\mbox{ a.s., }
\end{equation}
and the allele frequency spectrum for large samples is given by the following a.s. convergence 
\begin{equation}
\lim_{n\tendinfty}\frac{A_n(k)}{n} 	=	\intgen dx \,\theta\,e^{-\theta x} \frac{1}{W_\theta(x)^2}\left(1-\frac{1}{W_\theta(x)}\right)^{k-1}.
\end{equation}
Before ending this last subsection, we want to point out that in some cases, more explicit formulae can be computed. First, for the \emph{Yule process} with birth rate 1, (or with parameter $a$, but after replacing $\theta$ with $a\theta$), that is, when $W(x) = e^x$, one gets easily
$$
\lim_{n\tendinfty}\frac{S_n}{n}=\theta\quad\mbox{ and }\quad \lim_{n\tendinfty}\frac{S_n(k)}{n}=\frac{\theta}{k(k+1)}\ .
$$
Computations are not as straightforward for the number of haplotypes.
Second, for the \emph{critical birth--death process} with birth rate 1 (or with parameter $a$, but after replacing $\theta$ with $a\theta$),  that is, when $W(x) = 1+x$, one gets
$$
\lim_{n\tendinfty}\frac{S_n}{n \ln(n)}=\theta\quad\mbox{ and }\quad \lim_{n\tendinfty}\frac{S_n(k)}{n}=\frac{\theta}{k}\ .
$$
In addition,
$$
\lim_{n\tendinfty}\frac{A_n}{n}=\theta \ln\left(1+\theta^{-1}\right)\quad\mbox{ and }\quad \lim_{n\tendinfty}\frac{A_n(k)}{n}=\frac{\theta}{k} \left(1+\theta\right)^{-k}\ .
$$
\begin{rem}
It is amusing to notice that the rescaled number $A_n(k)$ of haplotypes with $k$ representatives is also the probability that a species has $k$ representatives in Fisher's \emph{log-series of species abundance} \cite{FCW}. In Fisher's model, a given species has an unknown density which is assumed to be drawn from a Gamma distribution with parameter $a$. As a result of Bernoulli sampling in a large population, it is then assumed that given the value $d$ of this density, the number $X$ of individuals spotted from this species is Poisson with parameter $\rho d$, where $\rho$ is the sampling intensity.
It can then be shown that as $a\downarrow 0$, conditional on $\{X\ge 1\}$ (since at least one individual must be spotted for the species to be recorded), $\PP(X=k)$ goes to $C(1+1/\rho)^{-k}/k$, for some normalising constant $C$.
\end{rem}

\begin{rem}
In a coalescent point process, divergence times are on average deeper than in the Kingman coalescent (our trees are more `star-like'). This forbids convergence of our statistics without rescaling (by the sample size $n$ or by $n\ln(n)$). In particular, notice that the asymptotic proportion of individuals in a cluster of size greater than $K$, i.e. $\lim_n n^{-1}\sum_{k\ge K}A_n(k)$, vanishes as $K$ grows to $\infty$. This shows that the largest cluster in a sample of $n$ has neglectable size w.r.t. $n$, which contrasts with the Kingman coalescent, where the allele frequency spectrum is given by Ewens' sampling formula (see \cite{D1, Ew}). As $n\to\infty$, the numbers of haplotypes $A_n(k)$ carried by $k$ individuals \cite{ABT} converge to independent Poisson r.v. with parameter $\theta/k$, and the $i$-th eldest haplotype \cite{DT} is carried by approximately $P_in$ individuals, where $(P_i;i\ge 1)$ is a Poisson--Dirichlet r.v. 
\end{rem}

\section{Number of polymorphic sites}
\label{sec : S}

Results for polymorphic sites depend on integrability assumptions on $H$. Of course these are always fulfilled if the time $t$ when the population was founded is known, since then $H\le t$ a.s. 
We will see that the critical assumptions are either $\EE(\min(H_1,H_2))<\infty$, or the more stringent $\EE (H)<\infty$.
Notice that the first assumption is equivalent to the integrability of $1/W^2$, and the second one to the integrability of $1/W$.

\subsection{Law of large numbers and central limit theorem}

Recall that $S_n$ is the number of polymorphic sites in the sample of $n$ individuals.

\begin{thm}
\label{thm : LLN+CLT}
If $\EE(H)<\infty$, then
$$
\lim_{n\tendinfty} n^{-1}S_n = \theta\, \EE(H)\qquad\mbox{a.s. and in }L^1.
$$
If in addition $\EE(H^{2})<\infty$, then 
$$
\sqrt{n}\left(n^{-1}S_n-\theta\,\EE(H)\right)
$$
converges in distribution to a centered normal variable with variance $\theta\,\EE(H)+\theta^2\mbox{Var}(H)$.
\end{thm}

\paragraph{Proof.}
Set $Y_n:=\max\{H_1, \ldots, H_{n-1}\}$.
Recall from the Introduction that 
$$
S_n=\sum_{i=1}^{n-1} Q_i +R_n,
$$
where $Q_i$ is the number of points of the Poisson point process ${\cal P}_i$ in $(0,H_i)$, and $R_n$ is the number of points of the Poisson point process ${\cal P}_0$ in $(0,Y_n)$.
By the strong law of large numbers, we know that
$$
\lim_{n\tendinfty} n^{-1}\sum_{i=1}^{n-1} Q_i = \theta\, \EE(H)\qquad\mbox{a.s. and in }L^1,
$$
so we need to prove that
$$
\lim_{n\tendinfty} n^{-1} R_n = 0\qquad\mbox{a.s. and in }L^1.
$$
Now because $R_n/Y_n$ converges to $\theta$ a.s. and in $L^1$, it is sufficient to prove that
$$
\lim_{n\tendinfty} n^{-1} Y_n = 0\qquad\mbox{a.s. and in }L^1.
$$
Because $Y_n<\sum_{i=1}^{n-1}H_i$, 
$$
\limsup_{n\tendinfty} n^{-1} Y_n =:Y<\infty\mbox{ a.s.}
$$
By the 0-1 law, $Y$ is not random. To prove that $Y=0$, we let $Y_n^{(1)}$ (resp. $Y_n^{(2)}$) be the maximum of the $H_i$'s indexed by odd (resp. even) numbers. Then it is clear that $Y_n=\max(Y_n^{(1)}, Y_n^{(2)})$, and that $n^{-1}Y_n^{(1)}$ as well as $n^{-1}Y_n^{(2)}$ both converge to $Y/2$. This shows that $Y=Y/2$, so that $Y=0$.

For convergence in $L^1$, pick any $x>0$, and notice that
\debeq
n^{-1}\EE(Y_n)	&=&		n^{-1}\EE(Y_n,Y_n\le x)+n^{-1}\EE(Y_n, Y_n>x)\\
								&\le&	n^{-1}x+n^{-1}\EE\left(\sum_{i=1}^{n-1}H_i\indic{H_i>x}\right)\\
								&\le&	n^{-1}x+\EE(H,H>x).
\fineq
Since $\EE(H)<\infty$, this last inequality shows that $n^{-1}\EE(Y_n)$ vanishes as $n\tendinfty$.

Now we prove the central limit theorem for $S_n$. It is elementary to compute $\mbox{Var}(Q_1)$ as $\theta\,\EE(H)+\theta^2\mbox{Var}(H)$, so by the classical central limit theorem applied to the sum of $Q_i$'s, we only have to prove that $R_n/\sqrt{n}$ converges to 0 in probability. For any $\lbd>0$,
$$
\EE\left(\exp\left(-\lbd R_n/\sqrt{n}\right)\right)=\EE\left(\exp\left(-\theta Y_n\left(1-e^{-\lbd/\sqrt{n}}\right)\right)\right),
$$
which shows it is sufficient to prove that $Y_n/\sqrt{n}$ converges to 0 in probability. As previously, we write
\debeq
n^{-1}\EE\left(Y_n^2\right)	&=&		n^{-1}\EE\left(Y_n^2,Y_n\le x\right)+n^{-1}\EE\left(Y_n^2, Y_n>x\right)\\
								&\le&	n^{-1}x^2+n^{-1}\EE\left(\sum_{i=1}^{n-1}H_i^2\indic{H_i>x}\right)\\	
								&\le&	n^{-1}x^2+\EE\left(H^2,H>x\right).
\fineq
Thus, convergence of $Y_n/\sqrt{n}$ to 0 holds in $L^2$, and subsequently, it holds in probability.
\hfill$\Box$

\subsection{Explicit formulae for the expected frequency spectrum}

Recall that $S_n(k)$ denotes the number of mutant sites that are carried by exactly $k$ individuals in the sample of $n$ individuals (and since we only count polymorphic sites, $S_n(n)=0$).

\begin{thm} 
\label{thm : exact exp}
For all $1\le k\le n-1$,
\begin{equation*}
\EE(S_n(k)) = \theta\,\intgen dx \left(1-\frac{1}{W(x)}\right)^{k-1}
\left( \frac{n-k-1}{W(x)^2}+\frac{2}{W(x)} \right),
\end{equation*}
which is finite if and only if $\EE(H)<\infty$.
Then in particular,
$$
\lim_{n\tendinfty}n^{-1}\EE(S_n(k)) = \theta\,\intgen \frac{dx}{W(x)^2}\left(1-\frac{1}{W(x)}\right)^{k-1}.
$$
\end{thm}
\begin{rem}
Taking the sum over $k$ in the r.h.s. of the last equality of the theorem, one gets $\theta\,\EE(H)$, so that, thanks to the $L^1$ convergence in Theorem \ref{thm : LLN+CLT},
$$
\lim_{n\tendinfty}n^{-1}\sum_{k=1}^{n-1}\EE(S_n(k))=\theta\,\EE(H)=\sum_{k\ge 1}\lim_{n\tendinfty} n^{-1}\EE (S_n(k)). 
$$ 
\end{rem}
Before giving a proof of the previous theorem, we want to make a point that will also be useful in the next subsection. For any tree with point mutations, a mutation is carried by $k$ individuals if and only of it is in the part of the tree subtending $k$ leaves. Then in any given tree with edge lengths and Poisson point process of mutations (with rate $\theta$) independent of the genealogy (as in our situation), the expectation of the number of mutations carried by $k$ individuals is $\theta L_k$, where $L_k$ is the Lebesgue measure of the part of the tree subtending $k$ leaves (i.e., tips).
In our setting, we will call $L_k(n)$, for $k\le n-1$, the Lebesgue measure of the part of the tree subtending $k$ tips among individuals $\{0,1,\ldots, n-1\}$, so that
$$
\EE(S_n(k)) = \theta\, \EE(L_k(n))\qquad 1\le k\le n-1.
$$

\begin{rem}
The last equality along with more specific considerations given in the next subsection provide a less analytic and more transparent proof than the proof we give hereafter. However, we stick to it for the interest of the method itself.
\end{rem}

\paragraph{Proof of Theorem \ref{thm : exact exp}}
We set $N(x)$ to be the smallest $i\ge 1$ such that $H_i> x$. The proof relies on the fact that
$$
\EE(S_n(k)) = \lim_{x\tendinfty}\theta\,\EE(L_k(N(x)) \mid N(x)=n)\qquad  1\le k\le n-1.
$$
On the event $\{N(x)=n\}$, we will need to extend the definition of  $L_k(N(x))$ to $k=n$,  as being the Lebesgue measure of the part of the tree \emph{up to time $-x$} subtending all tips $\{0,1,\ldots, n-1\}$, that is, $L_n(N(x))= x-\max_{i=1,\ldots, n-1}H_i$.

For editing reasons, we will prefer to write $F(x)=\PP(H>x)$, instead of $1/W(x)$. Since $F$ is a.e. differentiable and our goal is to let $x\tendinfty$, we can set $f(x) := -F'(x)$ without loss of generality. 
We let $\tilde{H}$ denote the branch length $H_{N(x)}$, and we set
$$
\tilde{N}:=\min\{k\ge 1:H_{N(x)+k}>x+dx\}, 
$$ 
as well as $\tilde{L}_k$ the Lebesgue measure of the part of the tree subtending $k$ leaves among individuals $\{N(x),N(x)+1,\ldots, N(x+dx)-1\}$.
Note that $(\tilde{N}, \tilde{L}_k, \tilde{H})$ are independent of $(N(x), L_k(N(x)))$; that $\tilde{H}$ is distributed as $H$ conditional on $\{H>x\}$; and that ($\tilde{N},\tilde{L}_k)$ is independent of $\tilde{H}$ and distributed as $(N(x+dx), L_k(N(x+dx))$. Next observe that if $\tilde{H}>x+dx$, then $N(x+dx)=N(x)$ and $L_k(N(x+dx))=L_k(N(x))$, except if $k=n$, where by definition $L_n(N(x+dx))=L_n(N(x))+dx$. On the other hand, if $\tilde{H}\in dx$, $L_k(N(x+dx))$ is the sum of measures of edges subtending $k$ tips in $\{0,1,\ldots,N(x)-1\}$ with measures of edges subtending $k$ tips in $\{N(x),\ldots,N(x)+\tilde{N}-1\}$. This reads
\begin{multline*}
L_k(N(x+dx))\indic{N(x+dx)=n} = \indic{\tilde{H}>x+dx}\indic{N(x)=n}\left(L_k(N(x))+dx\indicbis{k=n}\right)\\
	+\indic{\tilde{H}\le x+dx}\sum_{j=1}^{n-1}\indic{N(x)=j}\indic{\tilde{N}=n-j}\left(L_k(N(x))+\tilde{L}_k(\tilde{N})\right),
\end{multline*}
where we have used the extension of the definition of $L_k$ specified earlier (cases when $k=j$ or $k=n-j$ in the sum). Now set
$$
U_{k,n}(x):= \EE(L_k(N(x)), N(x)=n).
$$
By the independences stated previously, taking expectations, we get
$$
U_{k,n}'(x+)= -U_{k,n}(x)\frac{f}{F}(x)+\indicbis{k=n}\PP(N(x)=k)
	+2\sum_{j=1}^{n-1}U_{k,j}(x)\PP(N(x)=n-j)\frac{f}{F}(x).
$$
Setting
$$
V_k(x;s):=\sum_{n\ge k}U_{k,n}(x)s^n\qquad s\in[0,1),
$$
and observing that $\vert U_{k,n}(x)\vert\le nx$, and (so) that $\vert U_{k,n}'(x)\vert \le c(x)n^2$ for some positive $c(x)$ independent of $k$ and $n$,  we get 
$$
\frac{\partial V_k}{\partial x}(x;s)= -\frac{f}{F}(x)V_k(x;s)+\PP(N(x)=k) s^k + 2 \frac{f}{F}(x)\sum_{n\ge k}s^n\sum_{j=1}^{n-1}U_{k,j}(x)\PP(N(x)=n-j).
$$
Since $U_{k,j}(x)=0$ when $j\le k-1$, the last term equals
\debeq
 2 \frac{f}{F}(x)\sum_{n\ge k+1}s^n\sum_{j=k}^{n-1}U_{k,j}(x)\PP(N(x)=n-j)&=& 2 \frac{f}{F}(x)\sum_{j\ge k}U_{k,j}(x)s^j\sum_{n\ge j+1}s^{n-j}\PP(N(x)=n-j)\\
 	&=& 2 \frac{f}{F}(x)V_k(x;s)\sum_{n\ge 1}s^{n}\PP(N(x)=n).
\fineq
As a consequence, we get the following differential equation
$$
\frac{\partial V_k}{\partial x}(x;s)= G_k(x;s)V_k(x;s)+\PP(N(x)=k)s^k,
$$
where we have put
$$
G_k(x;s):=\left(2\EE\left(s^{N(x)}\right)-1\right)\frac{f}{F}(x).
$$
Now since $\PP(N(x)=k)=F(x)(1-F(x))^{k-1}$, we easily get 
$$
\int_0^x G_k(y;s) \, dy = \ln\left[\frac{F(x)}{\left( 1-s+sF(x)\right)^2}\right].
$$
This allows us to integrate the differential equation in $V_k(.;s)$ to finally arrive at
$$
V_k(x;s)=\frac{s^k F(x)}{(1-s+sF(x))^2}
\int_0^x(1-F(y))^{k-1}(1-s+sF(y))^2\, dy.
$$
With the shortcuts $u:=1-F(x)$ and $v:=1-F(y)$, and using the series expansion of $(1-us)^{-2}$, we get 
$$
V_k(x;s)=s^k (1-u)\int_0^x v^{k-1} (1-vs)^2 \sum_{j\ge 1} ju^{j-1}s^{j-1}\, dy.
$$
It is elementary algebra to compute the following equality
$$
(1-vs)^2 \sum_{j\ge 1} ju^{j-1}s^{j-1} = 1+\sum_{j\ge 1}s^j u^{j-2}\left(j(u-v)^2+u^2-v^2\right),
$$
which yields
$$
V_k(x;s)= s^k(1-u)\int_0^x v^{k-1}\, dy+(1-u)\sum_{j\ge 1}\int_0^x v^{k-1}s^{k+j} u^{j-2}\left(j(u-v)^2+u^2-v^2\right)\, dy.
$$
Identifying this entire series with the definition of $V_k$, we get for all $1\le k\le n-1$,
\begin{multline*}
U_{k,n}(x)= F(x)(1-F(x))^{n-k-2}\int_0^x (1-F(y))^{k-1}\times\\
\times\left((n-k)(F(y)-F(x))^2+(1-F(x))^2- (1-F(y))^2 \right)\, dy.
\end{multline*}
As a consequence, 
\begin{multline*}
\EE(L_k(N(x)) \mid N(x)=n)=(1-F(x))^{-k-1}\int_0^x (1-F(y))^{k-1}\times\\
\times\left((n-k)(F(y)-F(x))^2+(1-F(x))^2- (1-F(y))^2 \right)\, dy.
\end{multline*}
which, by Beppo Levi's theorem, converges, as $x\tendinfty$, to
$$
\theta^{-1} \EE(S_n(k))=\intgen (1-F(y))^{k-1}
\left((n-k)F(y)^2+1- (1-F(y))^2 \right)\, dy,
$$
and this finishes the proof.\hfill$\Box$

\subsection{Site frequency spectrum of large samples}
Here, we assume that $\EE(\min(H_1,H_2))<\infty$, that is, $1/W^2$ is integrable.

\begin{thm} 
\label{thm : SFS large}
For all $1\le k\le n-1$, the following convergence holds a.s. (and in $L^1$ as well if $\EE(H)<\infty$)
\begin{eqnarray*}
\lim_{n\tendinfty}n^{-1}S_n(k) &=&	\theta\, \EE\left(\left( \min\{H_1,H_{k+1}\}-\max\{H_2,\ldots,H_k\}\right)^+\right)\\
 	&=&	\theta\,\intgen \frac{dx}{W(x)^2}\left(1-\frac{1}{W(x)}\right)^{k-1}.
\end{eqnarray*}
\end{thm}
\paragraph{Proof.}
Reasoning similarly as in the previous subsection, we see that a point mutation occurring on branch $i$ is carried by $k$ individuals if and only if it is carried by individuals $i,i+1,\ldots, i+k-1$, and by no one else. This happens if and only if this mutation, corresponding to the atom $\ell_{ij}$, say, of ${\cal P}_i$, has 
$$
\max\{H_{i+1},\ldots, H_{i+k-1}\}<\ell_{ij}<H_i,
$$
for the mutation to be carried by individuals $i,i+1,\ldots, i+k-1$, along with
$$
\ell_{ij}<H_{i+k},
$$
for the mutation not to be carried by others. 
More formally, we set ${\cal F}$ the space of point processes on $(0,\infty)$, and ${F}_k$ the set of $(k+1)$-dimensional arrays with values in ${\cal F}\times (0,\infty)$. Next, for any $\Xi\in F_k$, written as $\Xi = ((p_0,x_0),\ldots, (p_k, x_k))$ we define
$$
G(\Xi):= \mbox{Card}\left(p_0\cap \left(\max\{x_{1},\ldots, x_{k-1}\}, \min\{x_0, x_{k}\}\right)\right),
$$
where it is understood that the interval $(a,b)$ is empty if $a\ge b$. Then the number of mutations carried by $k$ individuals among the first $n$ can be written as
$$
S_n(k)=\sum_{i=0}^{n-k}G(\Xi_i),
$$
where 
$$
\Xi_i:=(({\cal P}_i, H_i),\ldots,({\cal P}_{i+k}, H_{i+k}))
$$
and, for the last term of the sum to be correctly written, $H_{n}$ is set to $+\infty$ (as $H_0$). Next, observe that 
$$
\EE(G(\Xi_1))	=\theta\, \EE\left(\left( \min\{H_1,H_{k+1}\}-\max\{H_2,\ldots,H_k\}\right)^+\right),
$$
so that $G(\Xi_1)$ is integrable (assumption stated before the theorem). Now for any $0\le r\le k$, the random values $G(\Xi_i)$, for $i$ such that $i=r\; [k+1]$ (standing for mod $(k+1)$), are i.i.d. and integrable, so by the strong law of large numbers, we have the following a.s. convergence 
$$
\lim_{n\tendinfty} n^{-1}\sum_{0\le i= r [k+1]\le n-k}G(\Xi_i)=\frac{1}{k+1}\;\EE(G(\Xi_1)).
$$ 
Actually, the convergence would also hold in $L^1$ if we had discarded mutations carried by individual $0$ and individual $n-k$, which involve terms that are not integrable if $\EE(H)=\infty$. If $\EE(H)<\infty$, then convergence holds in $L^1$.
%
%
Summing over $r$ these $k+1$ equalities, we get the convergence of $n^{-1}S_n(k)$ to $\EE(G(\Xi_1))$, and 
\debeq
\EE(G(\Xi_1))	&=&\theta\, \EE\left(\left( \min\{H_1,H_{k+1}\}-\max\{H_2,\ldots,H_k\}\right)^+\right)\\
							&=&\theta\,\EE\intgen dx\,\indicbis{x<\min\{H_1,H_{k+1}\}}\,\indicbis{x>\max\{H_2,\ldots,H_k\}}\\
							&=&\theta\,\intgen dx\,\PP(H>x)^2\,\PP(H<x)^{k-1},
\fineq
which ends the proof.\hfill$\Box$

\subsection{Stable laws}

Here, we tackle the case when $H$ is in the domain of attraction of a stable law, which happens in particular for a splitting tree  whose contour process is a stable Lévy process with no negative jumps with index $\alpha\in(1,2]$. If such a population is censused with intensity $c>0$ then the corresponding function $W$ (see Introduction) is
$$
W(x)=1+cx^{\alpha-1}\qquad x\ge 0.
$$
From now on, we will assume that $W$ has the form given in the foregoing display.
Recall that $1/W(x)$ is the probability that a branch has length greater than $x$. 
Observe that here $H$ is not integrable, so that Theorems \ref{thm : LLN+CLT} and \ref{thm : exact exp} do not apply. However, asymptotic results for the site frequency spectrum of large samples given in Theorem \ref{thm : SFS large} apply for $\alpha>3/2$.

\subsubsection{Brownian case}
Here, we assume that $\alpha=2$, which  corresponds both to a (censused) Brownian population and to the (censused or not) population of a critical birth--death process.
\begin{thm}
\label{thm : brown}
When $W(x)=1+cx$, we have the following convergence in probability
$$
\lim_{n\tendinfty} \frac{S_n}{n\ln (n)} = \theta/c.
$$
\end{thm}
\paragraph{Proof.}
Recall that $S_n$ is to be written as
$$
S_n=\sum_{i=1}^{n-1} Q_i +R_n,
$$
where $Q_i$ is the number of points of the Poisson point process ${\cal P}_i$ in $(0,H_i)$, and $R_n$ is the number of points of the Poisson point process ${\cal P}_0$ in $(0,Y_n)$, where $Y_n=\max\{H_1, \ldots, H_{n-1}\}$.
Now observe that
$$
\PP(Y_n>\vareps n\ln (n)) = 1-\left(1-\frac{1}{1+c\vareps n\ln (n)} \right)^{n-1},
$$
which vanishes as $n\tendinfty$, so that $Y_n/n\ln (n)$ converges to 0 in probability. This implies in turn that $R_n/n\ln (n)$ also converges to 0 in probability. As a consequence, we can focus on the sum of $Q_i$'s. Pick any $\lbd>0$ and check that
$$
\EE\left(\exp-\frac{\lbd}{n\ln(n)} \sum_{i=1}^{n-1} Q_i\right)=\left(\EE\left(\exp-\theta H\left(1-e^{-\lbd/n\ln(n)}\right)\right)\right)^{n-1},
$$
We are bound to study the behaviour of $\EE(\exp-yH)$ as $y\to 0$.
\debeq
\EE(\exp-yH)	&=& 1-y\int_0^\infty\frac{e^{-yx}}{W(x)}\, dx\\
							&=&	1-y\int_0^\infty\frac{e^{-u}}{y+cu}\, du\\
							&=&	1-y\int_1^\infty\frac{e^{-u}}{y+cu}\,du +y \int_0^1\frac{1-e^{-u}}{y+cu}\,du- yc^{-1}\ln((y+c)/y)\\
							&=& 1+c^{-1}y\ln(y) + O(y),
\fineq
where $O(y)/y$ is bounded near 0. Setting $u_n:=\theta\left(1-e^{-\lbd/n\ln(n)}\right)$, there is a vanishing sequence $v_n$ such that
\debeq
\EE\left(\exp-\lbd\frac{S_n}{n\ln(n)} \right) 		&=&\left(1+c^{-1}u_n\ln(u_n)+O(u_n)\right)^n(1+v_n)\\
	&=&\exp\left(c^{-1}nu_n\ln(u_n)+O(nu_n)\right)(1+v_n),		
\fineq
which converges to $\exp(-\lbd\theta/c)$.\hfill$\Box$

\subsubsection{Stable case $\alpha\not= 2$}
Here, we assume that $W(x)=1+cx^{\alpha-1}$, for some $\alpha\in(1,2)$.
\begin{thm}
\label{thm : stable}
When $W(x)=1+cx^{\alpha-1}$, we have the following convergence in distribution
$$
\lim_{n\tendinfty} \frac{S_n}{n^{1/(\alpha-1)}} = Z_{\varphi(\mathbf{e})},
$$
where $(Z_t;t\ge 0)$ is the stable subordinator with Laplace exponent $\lbd\mapsto c^{-1}\theta^{\alpha-1}\lbd^{\alpha - 1}$, $\mathbf{e}$ is an independent exponential r.v. with parameter 1, and $\varphi$ is defined by
$$
\varphi(x)=x^{1-\alpha}\,e^{-x}+\int_0^x ds\, s^{1-\alpha}\,e^{-s}\qquad x>0.
$$
\end{thm}
\begin{rem}
Observe that $\varphi$ decreases on $(0,\infty)$ from $+\infty$ to a positive limit, equal to $\Gamma(2-\alpha)$. Also, recall that $S_n=\sum_{i=1}^{n-1} Q_i +R_n$, where $R_n$ is the extra contribution from the maximum branch length. Then it is possible to see by the same kind of proof as that of the theorem, that $\sum_{i=1}^{n-1} Q_i$ converges in distribution to $Z_{\Gamma(2-\alpha)}$. This indicates that, opposite to the Brownian case, the (double) contribution of the maximum branch length is not negligible here. 
\end{rem}
\paragraph{Proof.}
Let us compute the limiting distribution of $n^{-1/(\alpha -1)}(Y_n+\sum_{i=1}^{n-1}H_i)$, where $Y_n=\max\{H_1, \ldots, H_{n-1}\}$. Set $\beta:=1/(\alpha-1)$, as well as
$$
I_n(\lbd):=\EE\left(\exp-\lbd\, n^{-\beta}\left(Y_n+\sum_{i=1}^{n-1}H_i\right)\right).
$$
Then
$$
I_n(\lbd)=\intgen\PP(Y_n \in dz) e^{-2\lbd n^{-\beta} z} \left( \EE\left(e^{-\lbd n^{-\beta}H_z'}\right)\right)^{n-2},
$$
where $H_z'$ has the law of $H$ conditioned on being smaller than $z$. Next, we have
$$
\PP(Y_n\in dz) = \left(\frac{cz^{\alpha-1}}{1+cz^{\alpha-1}}\right)^{n-2}\frac{c(n-1)(\alpha-1)z^{\alpha-2}}{(1+cz^{\alpha-1})^2}\,dz\qquad z>0
$$
and
$$
\PP(H_z'\in dx) = \frac{c(\alpha-1)x^{\alpha-2}}{(1+cx^{\alpha-1})^2}\,\frac{1+cz^{\alpha-1}}{cz^{\alpha-1}}\,dx\qquad 0< x < z,
$$
so we get
$$
I_n(\lbd)=\intgen dz \,\frac{c(n-1)(\alpha-1)z^{\alpha-2}}{(1+cz^{\alpha-1})^2}e^{-2\lbd n^{-\beta} z} \left( \int_0^z dx\, \frac{c(\alpha-1)x^{\alpha-2}}{(1+cx^{\alpha-1})^2}\,e^{-\lbd n^{-\beta}x}\right)^{n-2}
$$
Changing variables, this also reads
$$
I_n(\lbd)=c^{-1} (1-n^{-1})(\alpha-1)\lbd^{\alpha-1}\intgen dv \,\frac{v^{-\alpha} \,e^{-2v}}{\left(1+n^{-1}c^{-1}\lbd^{\alpha-1}v^{1-\alpha}\right)^2} J_n(v;\lbd)^{n-2}
$$
where 
\debeq
J_n(v;\lbd)	&=& (\alpha-1)cn\lbd^{1-\alpha}\int_0^v du\, \frac{u^{\alpha-2}\, e^{-u}}{(1+cn\lbd^{1-\alpha}u^{\alpha-1})^2}\\
						&=&	\left[\frac{- e^{-u}}{1+cn\lbd^{1-\alpha}u^{\alpha-1}}\right]_0^v- \int_0^v du\,\frac{ e^{-u}}{1+cn\lbd^{1-\alpha}u^{\alpha-1}}\\
						&=&1-\frac{ e^{-v}}{1+cn\lbd^{1-\alpha}v^{\alpha-1}}- \int_0^v du\,\frac{ e^{-u}}{1+cn\lbd^{1-\alpha}u^{\alpha-1}}\\
						&=& 1-n^{-1}K_n(v;\lbd),
\fineq
where $K_n(v;\lbd)$ is positive and converges to $c^{-1}\lbd^{\alpha-1}\varphi(v)$ as $n\tendinfty$.
By the Lebesgue convergence theorem, we get the convergence of $I_n(\lbd)$ to
$$
 c^{-1}(\alpha-1)\lbd^{\alpha-1}\intgen dv \,v^{-\alpha} \,e^{-2v} \exp (-c^{-1}\lbd^{\alpha-1}\varphi(v)).
$$
Integrating by parts with $\varphi'(v)=(1-\alpha)v^{-\alpha}e^{-v}$, we finally get
$$
\lim_{n\tendinfty} I_n(\lbd) = \intgen dv  \,e^{-v} \exp (-c^{-1}\lbd^{\alpha-1}\varphi(v)).
$$
The last step is the same as in the foregoing proof, that is
\debeq
\lim_{n\tendinfty} \EE\left(\exp-\lbd\, n^{-1/(\alpha-1)}S_n\right)
	&=&\lim_{n\tendinfty} \EE\left(\exp-\theta\,\left(1-e^{-\lbd n^{-1/(\alpha-1)}}\right)\left(Y_n+\sum_{i=1}^{n-1}H_i\right)\right)\\
	&=&\lim_{n\tendinfty} I_n(\theta\lbd )\\
		&=&\intgen dv  \,e^{-v} \exp (-c^{-1}\theta^{\alpha-1}\lbd^{\alpha-1}\varphi(v)),
\fineq
which is the desired result. \hfill $\Box$

\section{Number of distinct haplotypes}
\label{sec : A}

\subsection{The next branch with no extra mutation}
We let ${\cal E}^\theta$ denote the set of individuals who \emph{carry no more mutations than} individual $0$ (some of and at most exactly the mutations carried by $0$, but no other mutation). Set $K^\theta_0:=0$ and for $i\ge 1$, define $K^\theta_i$ as the $i$-th individual in ${\cal E}^\theta$, and $H^\theta_i:=H_{K^\theta_i}$ the associated branch length.
We write $H^\theta$ in lieu of $H^\theta_1$ and we define the function $W_\theta$ by
$$
\PP(H^\theta>x)=\frac{1}{W_\theta(x)}\qquad x\ge 0.
$$
\begin{prop}
\label{prop : next branch with no}
The bivariate sequence $((K^\theta_i-K^\theta_{i-1},H^\theta_i);i\ge 1)$ is a sequence of i.i.d. random pairs.
The function $W_\theta$ is given by
$$
W_\theta(x)=1+\int_0^x W'(u)e^{-\theta u}\, du\qquad x\ge 0.
$$
\end{prop}

\begin{rem}
In the case when the coalescent process is derived from a splitting tree with lifespan measure $\Lambda$, the calculation of $W_\theta$ is straightforward. Indeed, it can  be seen in that case that the point process $(H^\theta_i;i\ge 1)$ is the coalescent point process of the splitting tree obtained from the initial splitting tree with mutations after throwing away all points above a mutation. But this new tree is again a splitting tree, since lifespans are i.i.d. and terminate either at death time or at the first point mutation, so the lifespan measure is now $\Lambda_\theta (dx)=e^{-\theta x}\,\Lambda (dx)+ \theta e^{-\theta x} \Lambda((x,\infty))\, dx$. As a consequence, $W_\theta$ is here the scale function characterised as in \eqref{eqn : LT scale} by its Laplace transform
\debeq
\intgen dx\, e^{-\lbd x} \, W_\theta(x) &=& \left(\lbd -\intgen \Lambda_\theta(dx) (1-e^{-\lbd x}) \right)^{-1}\\
																				&=&\frac{\lbd+\theta}{\lbd}\left(\lbd+\theta -\intgen \Lambda(dx) (1-e^{-(\lbd +\theta) x}) \right)^{-1}\\
																				&=&\frac{\lbd+\theta}{\lbd}
\intgen dx\, e^{-(\lbd+\theta) x} \, W(x),
\fineq
which yields the equality given in the statement.
\end{rem}

\paragraph{Proof.}
First observe that the pair $(K^\theta_1, H^\theta_1)$ does not depend on the haplotype of individual $0$, and that the $i$-th individual with no mutation other than those carried by individual $0$ is also the next individual after $K^\theta_{i-1}$ with no mutation other than those carried by individual $K^\theta_{i-1}$. This ensures that $(K^\theta_i-K^\theta_{i-1},H^\theta_i)$ has the same law as $(K^\theta_1,H^\theta_1)$, and the independence between $(K^\theta_i-K^\theta_{i-1},H^\theta_i)$ and previous pairs is due to the independence of branch lengths and the fact that new mutations can only occur on branches with labels strictly greater than $K^\theta_{i-1}$.

Now the event $\{H^\theta\in dx\}$ can be decomposed according to: the value of $H_1$; conditional on $H_1=z$, the value of the age $V_z$ of the oldest mutation on $H_1$; conditional on $V_z=y$, the value $H_y'$ of the branch length associated with the first individual in ${\cal E}^\theta_1$ with branch length greater than $y$. Indeed, $H^\theta\in dx$ if: $H_1\in dx$ and there is no mutation in $H_1$ (then $K^\theta_0=1$); or $H_1\in dx$, the age of the oldest mutation on $H_1=x$ is $V_x=y<x$ and the next individual with no mutation other than those carried by individual $1$ and branch length $H_y'>y$ has $H_y'<x$; or $H_1=z<x$, the age of the oldest mutation on $H_1=z$ is $V_z=y<z$ and the next individual with no mutation other than those carried by individual $1$ and branch length $H_y'>y$ has $H_y'\in dx$.   
\begin{multline*}
\PP(H^\theta\in dx)= \PP(H_1\in dx)e^{-\theta x}+\PP(H_1\in dx)\int_0^x\PP(V_x\in dy)\PP(H_y'<x)\\
+\int_0^x\PP(H_1\in dz)\int_0^z\PP(V_z\in dy)\PP(H_y'\in dx). 
\end{multline*}
Thanks to the first statement of the proposition, $H_y'$ has the same law as $H^\theta$ conditioned on being greater than $y$. Then since $\PP(V_z\in dy)=\theta \,e^{-\theta(z-y)}\,dy$, we get
$$
\PP(H^\theta\in dx)
=\PP(H_1\in dx)(1-\PP(H^\theta>x) f(x))+\PP(H^\theta \in dx)\int_0^x \PP(H_1\in dz) f(z),
$$
where we have set
$$
f(x):=\int_0^x\, dy\, \theta \,e^{-\theta (x-y)}\,W_\theta (y)\qquad x\ge 0.
$$
We can drop the index 1 of $H_1$, since only its law now matters.
We can rewrite the last result as
$$
\PP(H\in dx)=\PP(H^\theta\in dx)(1-\int_0^x \PP(H\in dz) f(z))+\PP(H\in dx)\PP(H^\theta>x) f(x),
$$
which can be integrated as
$$
\PP(H > x)=\PP(H^\theta >x)(1-\int_0^x \PP(H\in dz) f(z)).
$$
Defining now the function $G$ as
$$
G(x):=\PP(H>x)(W_\theta(x)-f(x)),
$$
we get, thanks to the last integration,
$$
G(x)=1-\int_0^x \PP(H\in dz) f(z)-\PP(H>x)f(x).
$$
Integrating by parts yields
$$
G(x)=1-\int_0^x dz\,\PP(H>z) f'(z)= 1-\int_0^x dz\,\PP(H>z)(-\theta f(z) +\theta W_\theta (z))=1-\theta \int_0^x dz\, G(z), 
$$
which shows that $G(x)=e^{-\theta x}$. This reads
$$
W(x)=e^{\theta x} W_\theta(x)-\theta \int_0^x\, dy \,e^{\theta y}\,W_\theta (y). 
$$
One differentiation and one integration provide the result.\hfill $\Box$

\subsection{Main result}

\subsubsection{Statement}
Recall that $A_n(k)$ denotes the number of haplotypes carried by $k$ individuals in a sample of $n$.
\begin{thm}
\label{thm : AFS large}
For all $k\ge 1$, the following convergence holds a.s.
$$
\lim_{n\tendinfty} n^{-1}A_n(k) 	=	\intgen dx \,\theta\,e^{-\theta x} \frac{1}{W_\theta(x)^2}\left(1-\frac{1}{W_\theta(x)}\right)^{k-1}.
$$
In addition,
$$
\lim_{n\tendinfty} n^{-1}A_n 	=	\intgen dx \,\theta\,e^{-\theta x} \frac{1}{W_\theta(x)}=\EE\left(1-e^{-\theta H^\theta}\right).
$$
\end{thm}

Before proving this statement, we insert a (sub)subsection in which we state and prove a preliminary key result.

\subsubsection{The key lemma}

Recall that $\ell_{1i}$ denotes the (time elapsed since the) $i$-th (most recent) mutation on the first branch length. In particular, the mutations carried by individual 1 and not by individual 0 are exactly those $\ell_{1i}$ such that $\ell_{1i}<H_1$ (the other points of the process are thrown away).
Let $N_i$ denote the number of individuals whose \emph{most recent mutation} is $\ell_{1i}$.

\begin{lem}
\label{lem : key}
In an infinite sample, for any integer $k\ge 1$,
$$
\sum_{i\ge 1} \PP(N_i=k)= \intgen \theta \,e^{-\theta z} \, dz \,\frac{1}{W_\theta(z)^2}\left(1-\frac{1}{W_\theta(z)}\right)^{k-1}
$$
\end{lem}
\paragraph{Proof.}
In the first place, not to care for the fact that only mutations with $\ell_{1i}<H_1$ contribute, we consider the number $N_i'$ of individuals whose most recent mutation is $\ell_{0i}$, and we condition on $\ell_{0j}=v_j$, $j\ge 1$.
We will use later the fact that the law of $N_i$ conditional on $\ell_{1j}=v_j$, $j\ge 1$, is that of $N_i'\indicbis{v_i<H}$, where $H$ is independent of $N_i'$ and the point process $(\ell_{0i};i\ge 1)$.

Recall from the previous subsection that ${\cal E}^\theta$ is the set of individuals who carry no more mutations than individual 0, that $K_i^\theta$ is the $i$-th individual in ${\cal E}^\theta$, and $H_i^\theta:=H_{K_i^\theta}$. Then set $D_0:=0$ and
$$
D_i:= \inf\{j\ge 1 : H_j^\theta > v_{i-1}\}\qquad i\ge 1.
$$
Now observe that $N_i'=D_{i}-D_{i-1}$ for all $i\ge 1$ (for $N_1'$, the count includes individual 0). As an application of Proposition \ref{prop : next branch with no}, we get that conditional on $\ell_{0j}=v_j$, $j\ge 1$, 
$$
\PP(N_1'=k)= \PP(H^\theta <v_1 )^{k-1}\PP(H^\theta >v_1),
$$
whereas for any $i\ge 2$,
$$
\PP(N_i'\not=0)= \PP(H^\theta <v_i \mid H^\theta > v_{i-1})
\quad \mbox{ and }\quad \PP(N_i'=k \mid N_i' \not=0)= \PP(H^\theta <v_i )^{k-1}\PP(H^\theta >v_i).
$$
Recalling the relation between the laws of $N_i$ and $N_i'$ mentioned in the beginning of the proof, we get that conditional on $\ell_{1j}=v_j$, $j\ge 1$, 
$$
\PP(N_1=k)= \PP(H^\theta <v_1 )^{k-1}\PP(H^\theta >v_1)\PP(H >v_1).
$$
whereas for any $i\ge 2$,
$$
\PP(N_i\not=0)= \PP(H^\theta <v_i \mid H^\theta > v_{i-1})\PP(H>v_i).
$$
Now $\PP(N_i'=k \mid N_i' \not=0)=\PP(N_i=k \mid N_i \not=0)$, so we finally get (for $i\ge 2$)
\debeq
\PP(N_i=k)	&=&\PP(H^\theta <v_i )^{k-1}\PP(H^\theta <v_i \mid H^\theta > v_{i-1})\PP(H^\theta >v_i)\PP(H>v_i)\\	&=&\left(1-\frac{1}{W_\theta(v_i)}\right)^{k-1}\left(1-\frac{W_\theta(v_{i-1})}{W_\theta(v_i)}\right)\frac{1}{W(v_i)W_\theta(v_i)} .
\fineq
It is well-known that for the Poisson point process of mutations, 
$$
\PP(\ell_{1,i-1}\in dx, \ell_{1i} \in dz)= \frac{\theta^i x^{i-2}}{(i-2)!}\,e^{-\theta z}\, dx\, dz\qquad 0<x<z, i\ge 2,
$$
so that
\debeq
\sum_{i\ge 2}\PP(N_i=k)	&=& \sum_{i\ge 2} \intgen dz \int_0^z dx\,  \frac{\theta^i x^{i-2}}{(i-2)!} \, e^{-\theta z}\, 	\left(1-\frac{1}{W_\theta(z)}\right)^{k-1}\left(1-\frac{W_\theta(x)}{W_\theta(z)}\right)\frac{1}{W(z)W_\theta(z)}\\
												&=&\intgen dz\, \theta \, e^{-\theta z} 	\left(1-\frac{1}{W_\theta(z)}\right)^{k-1}\frac{1}{W(z)W_\theta(z)} \int_0^z dx\, \theta\,e^{\theta x} \left(1-\frac{W_\theta(x)}{W_\theta(z)}\right).	
\fineq
Now thanks to Proposition \ref{prop : next branch with no}, we can perform the following integration by parts on the last integral in the last display
\debeq
\int_0^z dx\, \theta\,e^{\theta x} \left(1-\frac{W_\theta(x)}{W_\theta(z)}\right)
				&=&	\left[ e^{\theta x} \left(1-\frac{W_\theta(x)}{W_\theta(z)}\right)\right]_0^z
+\frac{1}{W_\theta(z)}\int_0^z dx\, e^{\theta x} W_\theta'(x)\\
				&=&	-1+\frac{1}{W_\theta(z)}+	\frac{1}{W_\theta(z)}\int_0^z dx\, W'(x)\\
				&=&	\frac{W(z)}{W_\theta(z)}-1.
\fineq
This entails
$$
\sum_{i\ge 2}\PP(N_i=k) = \intgen dz\, \theta \, e^{-\theta z} 	\left(1-\frac{1}{W_\theta(z)}\right)^{k-1}\frac{1}{W(z)W_\theta(z)}\left(\frac{W(z)}{W_\theta(z)}-1\right). 
$$
But since 
$$
\PP(N_1=k)=\intgen dz\, \theta \, e^{-\theta z} 	\left(1-\frac{1}{W_\theta(z)}\right)^{k-1}\frac{1}{W(z)W_\theta(z)},
$$
the result follows.\hfill$\Box$

\subsubsection{Proof of Theorem \ref{thm : AFS large}}
For each individual $i\ge 0$, we denote by ${\cal A}_{ij}$ the set of individuals bearing the unique haplotype whose most recent mutation is $\ell_{ij}$. In particular, it is understood that ${\cal A}_{ij}=\emptyset$ whenever $\ell_{ij}>H_i$ (because no such haplotype exists). 

Now fix $M\ge 1$. Similarly as in the proof of Theorem \ref{thm : SFS large}, we can define
$$
G_M(\Xi_i):= \mbox{ Card } \{j\ge 1: \mbox{ Card } {\cal A}_{ij} \cap\{i,\ldots,i+M\} \ge k\}, 
$$
where
$$
\Xi_i:=(({\cal P}_i, H_i),\ldots,({\cal P}_{i+M}, H_{i+M})).
$$
Observe that $G_M$ is bounded from above, so that $G_M(\Xi_i)$ is integrable for all $i\ge 0$. 
Now for any $0\le r\le M$, the random variables $G_M(\Xi_i)$, for $i$ such that $i=r\; [M+1]$ (standing for mod $(M+1)$), are i.i.d. and integrable, so by the strong law of large numbers, we have the following convergence a.s. (and in $L^1$)
$$
\lim_{n\tendinfty} n^{-1}\sum_{0\le i= r [M+1]\le n-M}G_M(\Xi_i)=\frac{1}{M+1}\;\EE(G_M(\Xi_1)).
$$ 
Summing over $r$ these $M+1$ equalities, we get the following convergence a.s. (and in $L^1$)
$$
\lim_{n\tendinfty} n^{-1}\sum_{i=0}^{n-M}G_M(\Xi_i)=\EE(G_M(\Xi_1)).
$$ 
Our goal is now to let $M\tendinfty$. First define
$$
A_n'(k):=\sum_{i=0}^n\mbox{ Card } \{j\ge 1: \mbox{ Card } {\cal A}_{ij} \cap\{i,\ldots,n\} \ge k\}.
$$
Notice that 
$$
A_n'(k) = \sum_{h\ge k} A_n(h).
$$
Then for any $i=0,\ldots, n-M$, for any $j\ge 1$, if $\mbox{ Card } {\cal A}_{ij} \cap\{i,\ldots,i+M\} \ge k$, then $\mbox{ Card } {\cal A}_{ij} \cap\{i,\ldots,n\} \ge k$, so that $A_n'(k)\ge \sum_{i=1}^{n-M}G_M(\Xi_i)$, and
$$
\liminf_{n\tendinfty} n^{-1} A_n'(k)\ge \liminf_{n\tendinfty} n^{-1}\sum_{i=1}^{n-M}G_M(\Xi_i)= \EE(G_M(\Xi_1)).
$$
Letting $M\tendinfty$, Beppo Levi's theorem yields
$$
\liminf_{n\tendinfty} n^{-1} A_n'(k) \ge \EE\;\big[ \mbox{Card } \{j\ge 1: \mbox{ Card } {\cal A}_{1j} \ge k\}\big]= \sum_{j\ge 1} \PP(\mbox{Card } {\cal A}_{1j} \ge k)=:y_k.
$$
In the notation of the previous subsection $\mbox{Card } {\cal A}_{1j}=N_j$, so by Fubini--Tonelli's theorem, 
$$
y_k=\sum_{j\ge 1}\PP(N_j\ge k)=\sum_{j\ge 1}\sum_{h\ge k}\PP(N_j=h)=\sum_{h\ge k} x_h,
$$
where $x_k:= \sum_{j\ge 1}\PP(N_j=k)$. Thanks to Lemma \ref{lem : key} we have the following explicit expression for $x_k$
$$
x_k=\intgen \theta \,e^{-\theta z} \, dz \,\frac{1}{W_\theta(z)^2}\left(1-\frac{1}{W_\theta(z)}\right)^{k-1}.
$$
Now recall that $\sum_{k\ge 1}A_n'(k)=\sum_{h\ge 1}hA_n(h)=n$. Since it is easily seen that $\sum_{k\ge 1}y_k=\sum_{h\ge 1}hx_h=1$, by Fatou's lemma
$$
1=\sum_{k\ge 1}y_k \le \sum_{k\ge 1}\liminf_n n^{-1}A_n'(k)\le \liminf_n n^{-1}\sum_{k\ge 1} A_n'(k)=1.
$$
Then we would get a contradiction if there was $k_0$ such that $\liminf_n n^{-1}A_n'(k_0)>y_{k_0}$, so that for all $k\ge 1$ a.s.,
$$
\lim_{n\tendinfty} n^{-1}A_n'(k)=y_{k}.
$$
The first equation of the theorem stems from the fact that $A_n(k) = A_n'(k) -A_n'(k+1)$ and the second one by taking $k=1$ in the last display.
It takes an elementary integration by parts to check that
$$
y_1 =	\intgen dx \,\theta\,e^{-\theta x} \frac{1}{W_\theta(x)}=
\EE\left(1-e^{-\theta H^\theta}\right).
$$

\paragraph{Acknowledgments.} This work was partially funded by the project MAEV `Modèles Aléatoires de l'\'Evolution du Vivant'  of ANR (French national research agency).

\end{document}